\documentclass[10pt]{amsart}
\usepackage{color}

\definecolor{c20}{rgb}{0.,0.7,0.}
\definecolor{c30}{rgb}{0.,0.,1.}
\definecolor{c40}{rgb}{1,0.1,0.7}
\definecolor{c50}{rgb}{1,0,0}
\definecolor{c60}{rgb}{1,0.9,0.1}
\def\GH{\mathcal{F}}

\newcommand{\ve}{\varepsilon}

\newcommand{\abs}[1]{\left\lvert #1 \right\rvert}

\newcommand{\E}[1]{\mathbb{E}\left(#1\right)}

\newcommand{\pk}[1]{\mathbb{P} \left\{ #1 \right\} }

\newcommand{\R}{\mathbb{R}}

\newcommand{\BQN}{\begin{eqnarray}}
\newcommand{\EQN}{\end{eqnarray}}
\newcommand{\BQNY}{\begin{eqnarray*}}
\newcommand{\EQNY}{\end{eqnarray*}}

\newcommand{\BS}{\begin{sat}}
\newcommand{\ES}{\end{sat}}
\newcommand{\BT}{\begin{theo}}
\newcommand{\ET}{\end{theo}}
\newcommand{\BK}{\begin{korr}}
\newcommand{\EK}{\end{korr}}

\newcommand{\BD}{\begin{de}}
\newcommand{\ED}{\end{de}}

\newcommand{\BIT}{\begin{itemize}}
\newcommand{\EIT}{\end{itemize}}
\newcommand{\BDI}{\begin{description}}
\newcommand{\EDI}{\end{description}}

\newcommand{\BRM}{\begin{remarks}}
\newcommand{\ERM}{\end{remarks}}

\newcommand{\BEL}{\begin{lem}}
\newcommand{\EEL}{\end{lem}}

\newtheorem{theo}{Theorem}[section]
\newtheorem{sat}[theo]{Proposition}
\newtheorem{de}[theo]{Definition}
\newtheorem{lem}[theo]{Lemma}

\newtheorem{korr}[theo]{Corollary}
\newtheorem{remark}[theo]{Remark}
\newtheorem{remarks}[theo]{Remarks}

\newcommand{\nelem}[1]{{Lemma \ref{#1}}}

\newcommand{\netheo}[1]{{Theorem \ref{#1}}}

\newcommand{\prooftheo}[1]{ \textsc{\bf Proof of Theorem} \ref{#1}:}

\newcommand{\prooflem}[1]{\textsc{\bf Proof of Lemma} \ref{#1}:}
\newcommand{\proofkorr}[1]{\textsc{\bf Proof of Corollary} \ref{#1}:}

\newcommand{\COM}[1]{}

\newcommand{\QED}{\hfill $\Box$}

\def\rw{\rightarrow}

\def\IF{\infty}

\def\Cov{\mathrm{Cov}}

\date{}

\def\oo{(1+o(1))}

\def\LT{\left}
\def\RT{\right}
\def\H{\mathcal{H}}

\def\ooo{(1+o(1))}

\def\rw{\rightarrow}

\def\TT{\mathcal{T}}

\def\Del{\triangle}

\def\vn{\varepsilon}
\def\Var{\text{Var}}

\def\todis{\overset{d}\rightarrow}
\def\topro{\overset{p}\rightarrow}
\def\NN{\mathcal{N}}

\def\ovX{\overline X_H}

\def\ovSK{\overset{\leftarrow}K}
\def\ovSK{K^{\leftarrow}}
\def\ssG{self-similar Gaussian risk process\ }
\def\ssGs{self-similar Gaussian risk processes\ }

\usepackage{color}

\newcommand{\tb}[1]{{\textcolor{blue}{#1}}}

\newcommand{\Ea}[1]{{\textcolor{red}{#1}}}
\def\tb#1{#1}

\def\Ea#1{#1}
\newcommand{\lp}[1]{{\textcolor{blue}{#1}}}
\def\lp#1{#1}

\begin{document}

\title[Parisian ruin of s-s  Gaussian risk processes] {Parisian Ruin of Self-similar Gaussian risk processes}

\author{Krzysztof D\c{e}bicki}
\address{Krzysztof D\c{e}bicki, Mathematical Institute, University of Wroc\l aw, pl. Grunwaldzki 2/4, 50-384 Wroc\l aw, Poland}
\email{Krzysztof.Debicki@math.uni.wroc.pl}

\author{Enkelejd  Hashorva}
\address{Enkelejd Hashorva, University of Lausanne,\\
B\^{a}timent Extranef, UNIL-Dorigny, 1015 Lausanne, Switzerland
}
\email{Enkelejd.Hashorva@unil.ch}

\author{Lanpeng Ji }
\address{Lanpeng Ji, University of Lausanne,\\
B\^{a}timent Extranef, UNIL-Dorigny, 1015 Lausanne, Switzerland
}
\email{Lanpeng.Ji@unil.ch}

\bigskip

\date{\today}
 \maketitle

{\bf Abstract:} In this paper we derive the exact asymptotics of the probability of
Parisian ruin for self-similar Gaussian risk processes. \Ea{Additionally}, we obtain the normal
approximation of the Parisian ruin time and \tb{derive}
an asymptotic relation between the Parisian and the classical ruin times.



\section{Introduction}

Let $\{X_H(t), t\ge0\}$ be  a centered self-similar Gaussian process with almost surely continuous sample paths and index $H\in(0,1)$, i.e., $\Var(X_H(t))=t^{2H}$  and for any $a>0$ and $s,t\ge0$
\BQNY
\Cov(X_H(at),X_H(as))=a^{2H}\Cov(X_H(t),X_H(s)).
\EQNY
Let $\beta,c$ be two positive constants. In risk theory the surplus process of an insurance company can be modeled by
\BQN\label{eq:R}
R_u(t)=u+ct^\beta-X_H(t), \ \ \ t\ge0,
\EQN
where $u$ is the so-called initial reserve, $ct^\beta$ models the total premium received up to
time $t$, and $X_H(t)$ represents the total \tb{amount of aggregated
claims (including fluctuations) up to time $t$}. Typically, classical \Ea{risk} models \Ea{assume} a linear
premium income, meaning that $\beta=1$.
\tb{In} this paper we deal with \tb{a} more general case $\beta> H$ allowing for non-linear premium income. Below we shall refer to $R_u$  as the {\it self-similar Gaussian risk process}.
The justification for choosing self-similar processes to model the aggregated claim process comes from \cite{Mich98}, where it is
shown that the  ruin probability  for self-similar Gaussian risk processes  is a good approximation of the  ruin
probability for the classical risk process. 
 Recent contributions have shown that self-similar Gaussian processes such as fractional
 Brownian motion (fBm), sub-fractional Brownian motion and bi-fractional
 Brownian motion are useful in modeling of financial risks, see e.g.,
 \cite{EM02, HP99, HP08,KK04, RT06} and the references therein.

For any $u\ge0$, define the {\it classical ruin time} of the self-similar Gaussian risk process by
\BQN\label{eq:tau1}
\tau_u=\inf\{t\ge0: R_u(t)<0\} \ \ \ \lp{(\mathrm{with\ }\inf\{\emptyset\}=\IF)}
\EQN
and thus the {\it probability of ruin} is defined as
\BQN
\pk{\tau_u<\IF}.
\EQN
The classical ruin time and the probability of ruin for the self-similar Gaussian risk
process are well studied in the literature; see, e.g., \cite{HP99, HP08, dieker2005extremes}.

\tb{Recently,} an extension of the classical notion of ruin, that is \Ea{the}
{\it Parisian ruin}, focused substantial interest; see \cite{DW08, czarna2011ruin, Palm14a}
and \Ea{the} references therein. The core of the notion of \Ea{the} Parisian ruin is that
now one allows 
the surplus process to spend a pre-specified time under the level zero before
the ruin is recognized. To be more precise, let
$T_u$ model the pre-specified time which is a positive deterministic function of the initial
reserve $u$.  In our setup, the {\it Parisian ruin time} of  the \ssG $R_u$ is defined as
\BQN\label{eq:tau2}
\tau^*_u=\inf\{t\ge T_u: t-\kappa_{t,u}\ge T_u\},\ \ \ \text{with} \ \kappa_{t,u}=\sup\{s\in[0,t]: R_u(s)\ge0\}.
\EQN
Here we make the convention that $\sup\{\emptyset\}=0$.

In this contribution we focus on the Parisian ruin probability, i.e.,
\begin{eqnarray}
\pk{\tau^*_u<\IF}=\pk{\inf_{t\ge0} \sup_{s\in[t,t+T_u]} R_u(s)<0}.\label{par.ruin.pr}
\end{eqnarray}
We refer to \cite{czarna2011ruin,LCP13,Palmowski14, DW08}
for recent analysis of (\ref{par.ruin.pr}) for the L\'evy
surplus model.

Assume for the moment that $X_H$ is a standard Brownian motion, $\beta=1$ and $T_u= T>0, u>0$. Thus  $R_u$ is the  Brownian motion risk process with a linear trend. As shown in the deep contribution \cite{LCP13},
for any $u\ge0$
\BQN\label{deep}
\pk{\tau^*_u<\IF}= \frac{\exp\LT(- {c^2T}/{2}\RT)-c \sqrt{2 \pi T}\Phi(-c\sqrt{T})}{\exp\LT(-{c^2T}/{2}\RT)+c \sqrt{2 \pi T}\Phi( c\sqrt{T})} \exp(-2cu),
\EQN
where  $\Phi(\cdot)$ is the distribution of \tb{a standard Normal} random variable.
\tb{Since}
the case $\beta \neq1$ seems to be completely untractable, even for the Brownian motion
risk process, one has to resort \lp{to} bounds and asymptotic results,
allowing the initial capital $u$ to become large, see e.g., \cite{EKM97}.

\tb{This contribution is concerned with the asymptotic behaviour of the}
\lp{Parisian ruin probability} for a large class of \ssGs as $u\to\IF$.
Under a local stationary condition on the correlation of the self-similar
process $X_H$ (see \eqref{eq:locstaXH}) and a mild condition on
$T_u$ (see \eqref{eq:Tu}), in \netheo{Thm1} we derive the asymptotics of the
\lp{ Parisian ruin probability.} Interestingly, as a corollary, it appears that \lp{for the fBm risk process with a linear trend} if $H>1/2$, then
\BQN\label{main:A}
\pk{\tau^*_u<\IF}= \tb{\pk{\tau_u<\IF}(1+o(1))} , \quad u\to \IF
\EQN
even if $T_u$ grows to infinity at a specified rate, as $u\to\infty$.

\tb{The combination of the above with
the seminal contribution \cite{HP99}, where the exact asymptotic behaviour of $\pk{\tau_u<\IF}$
for a general class of \ssGs  was derived, \lp{implies thus} 
the exact asymptotic behaviour of
the \lp{Parisian ruin probability}.}

\tb{Additionally, we derive
the approximation for the conditional (scaled)
Parisian ruin time and the asymptotic relation between the classical ruin time
and the Parisian ruin time given that Parisian ruin occurs.
This results goes in line with,
e.g., \cite{AsmAlb10, DHJ13a, EKM97, HJ13, HP08, HZ08, GrifMal12, Grif13, HJ14c, Palm14b},
where the approximation of the classical ruin time is considered.
The obtained normal approximation of the Parisian ruin time
is a new result even for the Brownian motion risk process with a linear trend.
}

Brief outline of the paper:
\tb{
In Section 2 we introduce \Ea{our} notation. The exact asymptotics of \Ea{the} Parisian ruin probability
is given in Section \ref{s.par.ruin}; while the time of the Parisian ruin is analyzed
in Section \ref{s.ruin.time}.
Proofs of the above results are relegated to Section \ref{s.proofs}.
In Appendix we display \Ea{the} proof of Theorem \ref{ThmX}, that plays \lp{an} important role in the
proof of the main result.
}

\section{Notation}\label{s.notation}
Let \tb{ $\{X_H(t), t\ge 0\}$ be a \lp{centered} self-similar Gaussian process
with almost surely continuous sample paths and index $H\in(0,1)$,
as defined in the Introduction. By $\{B_\alpha(t), t\ge 0\}$ we denote a standard
 \lp{fBm}
with Hurst parameter $\alpha/2\in(0,1]$.}

\tb{It is useful to define, for $\beta>H$ and \Ea{$c>0$}}
\BQN\label{eq:X}
Z(t)=\frac{X_H(t)}{1+ct^\beta},\ \  t\ge0. 
\EQN
\tb{
Indeed, by self-similarity of $X_H$, for any $u$ positive
\BQN\label{eq:ss}
\pk{\tau^*_u<\IF}=
\pk{\sup_{t\ge0} \inf_{s\in[t,t+T_u]} \Bigl(X_H(s)-cs^\beta\Bigr)>u}=
\pk{\sup_{t\ge0}\inf_{s\in\LT[0,T_uu^{-\frac{1}{\beta}} \RT]} Z(t+s) >u^{1-\frac{H}{\beta}}}.
\EQN
Let
$\sigma_Z(t)= \sqrt{\Var(Z(t))}$ . It can be checked \lp{(cf. \cite{HP99, HP08})} that
$\sigma_Z(t)$
attains its maximum on $[0,\IF)$ at the unique point
$$t_0=\LT(\frac{H}{c(\beta-H)}\RT)^{\frac{1}{\beta}}$$
and
\[
\sigma_Z(t)=A-\frac{BA^2}{2}(t-t_0)^2 +o((t-t_0)^2)
\]
as $t\to t_0$, where
\BQN\label{eq:AB}
 A=\frac{\beta-H}{\beta} \LT(\frac{H}{c(\beta-H)}\RT)^{\frac{H}{\beta}},\ \ B=\LT(\frac{H}{c(\beta-H)}\RT)^{-\frac{H+2}{\beta}} H\beta.
\EQN
}

\tb{In the rest of the paper we assume} {\it the local stationarity}
of the standardized Gaussian process $\ovX(t):=X_H(t)/t^H, t\ge0$ in a
neighborhood of the point $t_0$ i.e.,
\BQN\label{eq:locstaXH}
\lim_{s\to t_0,t\to t_0}\frac{ \E{(\ovX(s)-\ovX(t))^2}}{K^2 (\abs{s-t})}=Q>0
\EQN
holds for some positive function $K(\cdot)$ which is regularly varying at 0 with index $\alpha/2\in(0,1)$.
Condition \eqref{eq:locstaXH} is common in the literature; most of the known self-similar Gaussian processes (such as fBm, sub-fBm, and bi-fBm) satisfy \eqref{eq:locstaXH}, see e.g., \cite{HJ13b}. Note that the local stationarity at $t_0$ and \Ea{the} self-similarity of the process $X_H$ imply the local stationarity of $X_H$
at any point $r>0$ i.e.,
$$
\lim_{s\to r, t\to r}\frac{ \E{(\ovX(s)-\ovX(t))^2}}{K^2 (\abs{s-t})}=\LT(\frac{t_0}{r}\RT)^\alpha Q.
$$
Throughout this paper we denote by $\ovSK(\cdot)$  the asymptotic inverse of  $K(\cdot)$; by definition 
 $$
 \ovSK(K(t))=K(\ovSK(t))\oo=t\oo,\ \ \ t\to0.
  $$
It follows that  $\ovSK(\cdot)$ is regularly varying at 0 with index ${2}/{\alpha}$; see, e.g., \cite{EKM97}.

\tb{
Let
$\H_{\alpha}$ be the classical Pickands  constant, \Ea{defined} by}
\[
\H_{\alpha}=\lim_{T\to\infty}\frac{1}{T}
\E{\exp \left(\sup_{t\in[0,T]}(\sqrt{2}B_\alpha(t)-t^\alpha)\right)}.
\]
\lp{We refer to} \cite{albin2010new, Berman92, debicki2002ruin, debicki2008note, demiro2003simulation,  DEJ14, DikerY, Man07, Pit96}
for the basic properties of {the} Pickands and related constants. \Ea{A new constant that shall appear in our results below is defined as}
\BQN\label{eq:GH}
\GH_{\alpha}(T)=\lim_{S\to\IF}\frac{1}{S}
\E{\exp\LT(\sup_{t\in[0,S]}\inf_{ s\in[0,T]}\LT(\sqrt{2}B_{\alpha}(t+s)-(t+s)^{\alpha}\RT)\RT)}\in(0,\IF)
\EQN
for any $T\in [0,\IF)$.

\section{Asymptotics of the Parisian ruin probability}\label{s.par.ruin}

\tb{In this section we display the main result of the paper, which is
the asymptotics of \Ea{the} Parisian ruin probability
$\pk{\tau^*_u<\IF}$, as $u\to\infty$,
for the self-similar Gaussian risk model introduced in Introduction.
}
\Ea{First, we note that} in the light of the seminal contribution \cite{HP99}
\BQN\label{eq:HPA1}
\pk{\tau_u<\IF}
&=&\frac{A^{\frac{3}{2}-\frac{2}{\alpha}} Q^{\frac{1}{\alpha}} \H_\alpha}{2^{ \frac{1}{\alpha}} B^{\frac{1}{2}}}\frac{ u^{  \frac{2H}{\beta}-2}}{ \overset{\leftarrow}K (u^{\frac{H}{\beta}-1}) }\exp\LT(-\frac{u^{2\LT(1-\frac{H}{\beta}\RT)}}{2A^2}\RT)\oo
\EQN
holds as $u\to\IF$. In order to control the growth of the deterministic time $T_u$,  \lp{ we shall} assume that
\BQN\label{eq:Tu}
\lim_{u\to\IF}\frac{T_u u^{-\frac 1\beta}}{\ovSK(u^{\frac{H}{\beta}-1})} =T\in[0,\IF).
\EQN

\BT\label{Thm1} Let $\{R_u(t), t\ge0\}$ be  the self-similar Gaussian risk process given as in \eqref{eq:R} with $X_H$ satisfying \eqref{eq:locstaXH} and $T_u, u>0$ satisfying \eqref{eq:Tu}. If  $\tau^*_u$ denotes the Parisian ruin time of $R_u$, then as $u\rw\IF$
\BQN\label{eq:main2}
\pk{\tau^*_u<\IF}= \frac{\GH_\alpha(D_0T)}{\H_\alpha} \pk{\tau_u<\IF}\ooo,
\EQN
where $D_0=2^{-\frac{1}{\alpha}} A^{-\frac{2}{\alpha}} Q^{\frac{1}{\alpha}}$ with $\GH_{\alpha}(T)$ defined in \eqref{eq:GH}.
\ET

The proof of Theorem \ref{Thm1} is deferred to Section \ref{s.proofs};
it relies on a general result for the asymptotics of
sup-inf functionals of Gaussian processes, given in \netheo{ThmX}.

\begin{remark}
\tb{
Observe that the Pickands constant $\H_\alpha=\GH_\alpha(0)$ and $\H_{1}=1$
(cf. \cite{Pit96}).
It is not clear how to calculate $\GH_{\alpha}(T)$
using the definition in \eqref{eq:GH}.
However for the special case $\alpha=1$,
\eqref{deep} and \eqref{eq:fbm} imply
}
\BQN\label{H12}
\GH _{1}( T)= 
\frac{\exp\LT(- {T}/{4}\RT)- \sqrt{ \pi T}\Phi(- \sqrt{T/2})}{\exp\LT(- {T}/{4}\RT)+ \sqrt{ \pi T}\Phi(- \sqrt{T/2})}
,\ \ T>0.
\EQN
In this paper we shall refer to  $\GH _{\alpha}(T)$ as the generalized Pickands constant.

\COM{For some constant  positive $v< D_0T$ we have
\footnote{K: maybe we can erase bound below?}
\BQNY
\GH _{1}(D_0T)&=& \lim_{S\to\IF}\frac{1}{S}
\E{\exp\LT(\sup_{t\in[0,S]}\inf_{ s\in[0,D_0T]}\LT(\sqrt{2}{B_{1}}(t+s)-(t+s)\RT)\RT)}\\
&\le& \lim_{S\to\IF}\frac{1}{S}
\E{\exp\LT(\sup_{t\in[0,S]}\LT(\sqrt{2}{B_{1/2}}(t+v)-\sqrt{2}{B_{1}}(v) -t +\sqrt{2}{B_{1}}(v)-v\RT)\RT)}\\
&=& \lim_{S\to\IF}\frac{1}{S}
\E{ \exp\LT(\sup_{t\in[0,S]}\LT(\sqrt{2}{B_{1}}(t) -t\RT)\RT)}
\E{ \exp\LT( \sqrt{2}\Bigl(B_{1}(v)-v\Bigr)\RT)}\\
&=&  1,
\EQNY
which in view of  \eqref{H12} is an obvious bound for $\GH_{1}(D_0T)$.
}
\end{remark}

As a  corollary of the last theorem we present next a result for the fBm  risk processes with a linear trend where $X_H$ is assumed to be a standard fBm. Specifically, we have
\BQNY
\Cov(X_H(t),X_H(s))=\frac{1}{2}(t^{2H}+s^{2H}-\mid t-s\mid^{2H}),\ \ \ t,s\ge0
\EQNY
and thus \eqref{eq:locstaXH} holds with $K(t)=t^H, t\ge0$ and $Q=t_0^{-2H} = [H/(c(\beta- H))]^{- 2H/\beta}$.

\BK\label{corofBm} Let
$R_u(t)=u+ct-B_{2H}(t), \lp{t\ge0}$
and let $T_u, u>0$ be such that
$\lim_{u\to\IF}T_u  u^{1/H-2}=T\in[0,\IF)$.
Then as $u\rw\IF$
\BQN\label{eq:fbm}
\pk{\tau^*_u<\IF}&=&  \GH_{2H}(D_0T) \frac{ 2^{-\frac{1}{2H}} }{\sqrt{H(1-H)}}\left(\frac{c^H u^{1-H}}{H^H (1-H)^{1-H}}\right)^{\frac{1}{H}-2}\nonumber\\
&&\times  \exp\left(-\frac{c^{2H} u^{2(1-H)}}{2H^{2H} (1-H)^{2(1-H)}}\right)(1+o(1)),
\EQN
where 
$D_0=2^{-\frac{1}{2H}} c^2H^{-2}  (1-H)^{2-\frac{1}{H}}$.
\EK

\begin{remark}
\tb{
Using \Ea{the fact} that $\GH_{2H}(0)=\H_{2H}$, Corollary \ref{corofBm}
implies that
\[\pk{\tau^*_u<\IF}=\pk{\tau_u<\IF}(1+o(1))\]
as $u\to\infty$,
if $T=0$ (i.e. $T_u=o(u^{(2H-1)/H})$).
Thus, if $H>1/2$, the asymptotics of \Ea{the} Parisian ruin probability
coincides with the asymptotics of the classical ruin
probability even if $T_u$  grows to infinity, \Ea{provided} that $T=0$.
This property is another manifestation of the long-range dependence structure of fBm with $H>1/2$.  
}\\
For the boundary case  $T_u= Tu^{ 1/H-2}$ with $T>0$, the Parisian ruin probability and the classical ruin probability
are not asymptotically equivalent, as the initial capital $u$ tends to infinity.
\end{remark}

In \cite{LRZ10} a different type of  Parisian ruin is considered, 
where the deterministic pre-specified time $T_u$ is
replaced by an independent random variable (in particular, an exponential random variable is dealt with therein, see also \cite{CzaPal14}).
In the following corollary we address the  Parisian ruin probability \lp{of this model.}

\BK\label{coro1} Let $\{R_u(t), t\ge0\}$ be  the self-similar Gaussian risk process given as in \eqref{eq:R} with $X_H$ satisfying \eqref{eq:locstaXH}. If  $\TT$ is a positive random variable independent of $\{R_u(t), t\ge 0\}$, then 
\BQN\label{eqcorr}
\pk{\inf_{t\ge0} \sup_{s\in[t,t+\TT]} R_u(s)<0}=
\pk{\tau_u<\IF}\ooo, \quad u\to \IF
\EQN
holds, provided that  $2H+\alpha>2\beta$.
\EK

\section{Asymptotics of \lp{the} Parisian ruin time}\label{s.ruin.time}

In this section we
present a normal approximation for the conditional (scaled)
Parisian ruin time. Additionally we derive an asymptotic
relation between the classical ruin time and the Parisian ruin time, given that the Parisian ruin occurs.\\
Hereafter $\todis$ and  $\topro$ stand for convergence in distribution and convergence in probability, respectively.

\BT\label{Thm2}
Let $\tau_u, \tau^*_u$ be the classical ruin time and the Parisian ruin time
\tb{for}
the self-similar Gaussian risk process $\{R_u(t), t\ge0\}$
given as in \eqref{eq:R}.
If  $X_H$ satisfies \eqref{eq:locstaXH} and $T_u, u>0$ satisfies \eqref{eq:Tu},  then as $u\rw\IF$
\BQN\label{eq:tausta}
\frac{\tau^*_u-t_0 u^{\frac{1}{\beta}}}{A^{\frac{1}{2}} B^{-\frac{1}{2}} u^{\frac{H}{\beta}+\frac1\beta-1}}\Big| (\tau^*_u<\IF) \todis \NN,
\EQN
where $A,B$ are as in \eqref{eq:AB} and $\NN$ is \tb{a standard Normal} random variable.
Moreover, as $u\to\IF$,
\BQN\label{eq:tautau}
\frac{\tau^*_u-\tau_u}{u^{\frac{H}{\beta}+\frac1\beta-1} }\Big| (\tau^*_u<\IF) \topro 0.
\EQN
\ET
The complete proof of Theorem \ref{Thm2} is given in Section \ref{s.proofs}.

As a straightforward implication of Theorem \ref{Thm2} \Ea{it follows that} if $H+1=\beta$, then
\BQN
(\tau^*_u-\tau_u) \Big| (\tau^*_u<\IF) \topro 0.
\EQN

\begin{remark}
In  \cite{HP08} a slightly more general class of Gaussian processes was considered.
Under additional technical conditions as A1 and A3 therein similar results as in \netheo{Thm1} and \netheo{Thm2} also hold for that class of Gaussian processes; the only difference is that in \eqref{eq:tausta} and \eqref{eq:tautau} we \lp{shall} have $\sqrt{\Var(X_H(u^{1/ \beta}))}$ instead of $u^{  H/ \beta }$ and $s_0(u)$ (in their notation) instead of $t_0$.\\
We note that extensions of our result  to Gaussian processes with random variance  under similar conditions as in \cite{MR2832415} are also possible.

\end{remark}

\section{Proofs}\label{s.proofs}
Before presenting the proofs of Theorems \ref{Thm1} and \ref{Thm2},
we shall derive a general result for the tail of sup-inf functional applied to
the Gaussian process $Z$.
Recall that in our notation
$\Phi(\cdot)$ is the distribution of an $N(0,1)$ random variable.

In what follows, \tb{in order to simplify the notation, we shall set}
\BQN\label{eq:qq}
q=q(v):=\ovSK\LT(\frac{1}{v}\RT),\ \ \ v>0.
\EQN

\BT\label{ThmX}
Let $\{Z(t), t\ge0\}$ be the centered Gaussian process given as in \eqref{eq:X}, and let $x_i(\cdot),i=1,2 $ be two functions such that
$\lim_{v\to\IF}x_i(v)=x_i, i=1,2$ and $\lim_{v\to\IF}x_i(v)v^{-1/2}=0, i=1,2$ for some $x_1, x_2\in \R\cup \{ \IF\}$ satisfying $x_2>-x_1$. Further, for all $v$ large set
$\Theta_{x_1,x_2}(v)=\LT[t_0-x_1(v) v^{-1}, t_0+x_2(v) v^{-1}\RT].$  
Then, for any positive function $\lambda(\cdot)$ such that $\lim_{v\to\IF}\lambda(v)=\lambda\in[0,\IF)$ we have, as $v\to\IF$
\BQN   \label{eq:I}
\pk{\sup_{t\in\Theta_{x_1,x_2}(v)}\inf_{s\in [0,\lambda(v) q ]}Z(t+s)>v}
&=&
\frac{\tb{\GH _\alpha(D_0\lambda)}}{\H _\alpha}
\LT(\Phi\LT(A^{-\frac{1}{2}} B^{\frac{1}{2}}x_2\RT)-\Phi\LT(-A^{-\frac{1}{2}} B^{\frac{1}{2}}x_1\RT)\RT)\nonumber\\
&&\times \pk{\sup_{t\ge0} Z(t)>v}\oo,
\EQN
where \tb{$D_0=2^{-\frac{1}{\alpha}} A^{-\frac{2}{\alpha}} Q^{\frac{1}{\alpha}}$,
and
$\GH _\alpha(\cdot)$}
defined in \eqref{eq:GH} \lp{is positive and finite}.
\ET
\tb{The complete proof of Theorem \ref{ThmX} is given in the Appendix.}

\tb{The next \Ea{result}
plays an important role in the proof of Theorem \ref{Thm1}.
We refer to \cite{HP99} for its proof.
}
\begin{lem}\label{l.pit}
Let $\{Z(t),t\ge0\}$ be defined as in \eqref{eq:X} and set $v(u)=u^{1- {H}/{\beta}}$. Then for any $G> t_0$ we have,
as $u\to\IF$
\BQN\label{eq:HPA2}
\quad \pk{\tau_u<\IF}&=&\pk{\sup_{t\in[0,G]}   \Bigl(X_H(t)-ct^\beta\Bigr)>u}\oo\nonumber\\
&=&
\pk{\sup_{t\in\LT[t_0-\frac{\ln v(u)}{v(u)},t_0+\frac{\ln v(u)}{v(u)}\RT]} Z(t)>v(u)}\oo.
\EQN
Further, as $u\to\IF$
\BQN\label{eq:HPA3}
 \pk{\sup_{\abs{t-t_0}>\frac{\ln v(u)}{v(u)}} Z(t)>v(u)}= o\LT(\pk{\sup_{t\ge0} Z(t)>v(u)} \RT).
\EQN
\end{lem}

\prooftheo{Thm1} 
The proof is based on an application of \netheo{ThmX}.
\tb{From \eqref{eq:ss} we straightforwardly have that}
\BQNY
\pk{\tau^*_u<\IF}=\pk{\sup_{t\ge0}\inf_{s\in\LT[0, S_v \RT]} Z(t+s)>v},
\EQNY
with
$$
 v=v(u)=u^{1-\frac{H}{\beta}}\ \ \ S_{v}=S_{v(u)}= {T_u}{u^{-\frac{1}{\beta}}}, \ \ u>0.
$$
\lp{Further}, condition \eqref{eq:Tu} implies $\lim_{v\to\IF} S_{v} /q=T\in[0,\IF)$, and  
 \BQN 
\Pi(v)\le \pk{\sup_{t\ge0}\inf_{s\in\LT[0, S_v \RT]} Z(t+s)>v} \le  \Pi(v)+\Sigma(v),
\EQN
where
\BQNY
\Pi(v)&=&\pk{\sup_{t\in\LT[t_0-\frac{\ln v}{v}, t_0+\frac{\ln v}{v}\RT]}\inf_{s\in\LT[0,S_v\RT]} Z(t+s)>v}\\
\Sigma(v)&=& \pk{\sup_{\abs{t-t_0}\ge \frac{\ln v}{v}} Z(t)>v}.
\EQNY
Taking $x_1(v)=x_2(v)=\ln v$ and $\lambda(v)=S_v /q$ in   \netheo{ThmX} we conclude  that as $u\to\IF$
$$
 \Pi(v)= \frac{\GH _\alpha(\lp{D_0}T)}{\H_\alpha} \pk{\sup_{t\ge0}  Z(t )>v}\oo=\frac{\GH _\alpha(\lp{D_0}T)}{\H_\alpha} \pk{\tau_u<\IF}\oo.
$$
Moreover, from \eqref{eq:HPA3} we have as $u\to\IF$
$$
\Sigma(v)=o(\Pi(v))
$$
establishing the proof.
 \QED

\proofkorr{coro1} For any $u>0$ we have 
\BQNY
\pk{\sup_{t\ge0} \inf_{s\in[t,t+\TT]} \Bigl( X_H(s)-cs^\beta\Bigr)>u}&\le& \pk{\sup_{t\ge0}  \Bigl( X_H(s)-cs^\beta \Bigr)> u}=\pk{\tau_u<\IF}.
\EQNY
Further, for any small positive $\vn\in(0, 2H+\alpha-2\beta) $ by the independence
of $\mathcal{T}$ and \tb{the} risk process
\BQNY
\lefteqn{\pk{\sup_{t\ge0} \inf_{s\in[t,t+\TT]} \Bigl(X_H(s)-cs^\beta\Bigr)>u}}\\
&\ge& \pk{\sup_{t\ge0}\inf_{s\in[t,t+\TT]} \Bigl(X_H(s)-cs^\beta\Bigr)>u, \TT<u^{\frac{2H+\alpha-2\beta-\vn}{\alpha\beta }}}\\
&\ge& \pk{\sup_{t\ge0}\inf_{s\in[t,t+u^{\frac{2H+\alpha-2\beta-\vn}{\alpha\beta }}]} \Bigl(X_H(s)-cs^\beta \Bigr)>u}\pk{ \TT<u^{\frac{2H+\alpha-2\beta-\vn}{\alpha\beta }}}.
\EQNY
Hence, the claim follows from \netheo{Thm1}, by letting $u\to\IF$. \QED

\prooftheo{Thm2} We use the same notation as in the proof of \netheo{Thm1}. For any $x\in\R$ and $u>0$
\BQNY
\pk{\tau^*_u<\IF}\pk{\frac{\tau^*_u-t_0 u^{\frac{1}{\beta}}}{A^{\frac{1}{2}} B^{-\frac{1}{2}} u^{\frac{H}{\beta}+\frac1\beta-1}}\le x\Big|\tau^*_u<\IF}&=&
\pk{\tau^*_u\le t_0 u^{\frac{1}{\beta}}+A^{\frac{1}{2}} B^{-\frac{1}{2}} x u^{\frac{H}{\beta}+\frac1\beta-1} }.
\EQNY
Next we focus on the asymptotics of
\BQNY
\pk{\tau^*_u\le t_0 u^{\frac{1}{\beta}}+A^{\frac{1}{2}} B^{-\frac{1}{2}} x u^{\frac{H}{\beta}+\frac1\beta-1} }
&=&\pk{\sup_{t\in[0,  t_0 u^{\frac{1}{\beta}}+A^{\frac{1}{2}} B^{-\frac{1}{2}} x u^{\frac{H}{\beta}+\frac1\beta-1} ]} \inf_{s\in[t,t+T_u]} \Bigl(X_H(s)-cs^\beta\Bigr)>u}\\
&=&\pk{\sup_{t\in[0,  t_0  +  A^{\frac{1}{2}} B^{-\frac{1}{2}} x  v^{-1} ]} \inf_{s\in[0,S_v]} Z(t+s)>v},
\EQNY
where 
$$
 v=v(u)=u^{1-\frac{H}{\beta}},\ \ \ S_{v}=S_{v(u)}= {T_u}{u^{-\frac{1}{\beta}}}, \ \ u>0.
$$
Similarly to the proof of \netheo{Thm1}, we have 
 \begin{eqnarray*} 
\Pi_0(v)\le \pk{\sup_{t\in[0,  t_0  +  A^{\frac{1}{2}} B^{-\frac{1}{2}} x  v^{-1} ]} \inf_{s\in[0,S_v]} Z(t+s)>v} \le  \Pi_0(v)+\Sigma_0(v),
\end{eqnarray*}
where
\BQNY
\Pi_0(v)&=&\pk{\sup_{t\in\LT[t_0-\frac{\ln v}{v}, t_0+ A^{\frac{1}{2}} B^{-\frac{1}{2}} x  v^{-1}\RT]}\inf_{s\in\LT[0,S_v\RT]} Z(t+s)>v}\\
\Sigma_0(v)&=& \pk{\sup_{ t\in[0, t_0-\frac{\ln v}{v}]  } Z(t)>v}.
\EQNY
In the light of \netheo{ThmX} and \eqref{eq:HPA3} we conclude that, as $u\to\IF$
\BQNY
\pk{\tau^*_u\le t_0 u^{\frac{1}{\beta}}+A^{\frac{1}{2}} B^{-\frac{1}{2}} x u^{\frac{H}{\beta}+\frac1\beta-1} }&=& \Ea{(1+o(1))}
\frac{\GH_\alpha (\lp{D_0}T)}{\H_\alpha}\pk{\tau_u<\IF}\Phi(x). 
\EQNY
Therefore, the claim of \eqref{eq:tausta} follows by applying \netheo{Thm1}. Moreover, as shown in  \cite{HP08}, Theorem 1
$$
\frac{\tau_u-t_0 u^{\frac{1}{\beta}}}{A^{\frac{1}{2}} B^{-\frac{1}{2}} u^{\frac{H}{\beta}+\frac1\beta-1}}\Big|(\tau_u<\IF)\todis \widetilde{\mathcal{N}},\ \ u\to\IF,
$$
with $\widetilde{\mathcal{N}}$  an $N(0,1)$
random variable. Consequently, by Lemma 2.3 in \cite{HJ13}
$$
\LT(\frac{\tau_u-t_0 u^{\frac{1}{\beta}}}{A^{\frac{1}{2}} B^{-\frac{1}{2}} u^{\frac{H}{\beta}+\frac1\beta-1}}, \frac{\tau^*_u-t_0 u^{\frac{1}{\beta}}}{A^{\frac{1}{2}} B^{-\frac{1}{2}} u^{\frac{H}{\beta}+\frac1\beta-1}}\RT)\Big|(\tau^*_u<\IF)\todis (\widetilde{\mathcal{N}},\widetilde{\mathcal{N}}),\ \ u\to\IF
$$
implying thus \Ea{\eqref{eq:tautau}}. This completes the proof.
\QED

\section{Appendix}

This Appendix is dedicated to the proof of \netheo{ThmX}.
We first present a crucial lemma which can be seen as an extension  of the celebrated Pickands lemma; see, e.g., \cite{PicandsA, Pit72, Pit96}. We refer to \cite{DebKo2013} for recent developments in this direction.

Let $\lambda_1,\lambda_2$ be two given positive constants. Consider the family of centered Gaussian random fields
$$ \{X_v(t,s), (t,s)\in[0,\lambda_1]\times [0,\lambda_2]\}$$
indexed by $v>0$. We shall assume that its variance equals 1 and the correlation functions
$r_v(t,s,t',s')=\Cov(X_v(t,s),X_v(t',s')), (t,s), (t',s')\in[0,\lambda_1]\times [0,\lambda_2], v>0$ satisfy the following two conditions:

{\bf C1}. There exist constants $D>0, \alpha\in(0,2]$ and  a positive function $f(\cdot)$ defined in $(0,\IF)$ such that
\BQNY
\lim_{v\to\IF}(f(v))^2(1-r_v(t,s,t',s'))=D\abs{s+t-s'-t'}^{\alpha}
\EQNY
holds for any $(t,s), (t',s')\in[0,\lambda_1]\times [0,\lambda_2]$;

{\bf C2}. There exist constants $C>0, v_0>0, \gamma\in(0,2]$ such that, for any $v>v_0$, with $f(\cdot)$ given in C1,
\BQNY
(f(v))^2(1-r_v(t,s,t',s'))\le C(\abs{s-s'}^\gamma+\abs{t-t'}^{\gamma})
\EQNY
holds uniformly with respect to  $(t,s), (t',s')\in[0,\lambda_1]\times [0,\lambda_2]$.

\BEL\label{LemGP} Let
$\{X_v(t,s), (t,s)\in[0,\lambda_1]\times [0,\lambda_2]\}, v>0$  be the family of
centered Gaussian random fields \tb{with variance  equal to 1 \Ea{defined above}}. If both {\bf C1} and {\bf C2} hold, then for any  positive function $\theta(\cdot)$
satisfying $\lim_{v\to\IF}f(v)/\theta(v)=1$ we have
\BQN\label{eq:GP}
\pk{\sup_{t\in[0,\lambda_1]}\inf_{ s\in[0,\lambda_2]} X_v(t,s)>\theta(v)}=\H_{\alpha}(D^{\frac{1}{\alpha}}\lambda_1,D^{\frac{1}{\alpha}}\lambda_2)
\frac{1}{\sqrt{2\pi}\theta(v)}\exp\LT(-\frac{(\theta(v))^2}{2}\RT)(1+o(1))
\EQN
as $u\to \IF$, where
\BQNY
\H_{\alpha}(\lambda_1,\lambda_2)=\E{\exp\LT(\sup_{t\in[0,\lambda_1]}\inf_{ s\in[0,\lambda_2]}\LT(\sqrt{2}B_{\alpha}(t+s)-(t+s)^{\alpha}\RT)\RT) } \in (0,\IF).
\EQNY
\EEL

\prooflem{LemGP} Note that the sup-inf functional satisfies {\bf F1-F2} in \cite{DebKo2013}. The proof follows by similar arguments as the proof of Lemma 1 therein, and therefore we omit the technical details. \QED

\def\lamm{\lambda^-_{\vn_0}}
\def\lamp{\lambda^+_{\vn_0}}

\prooftheo{ThmX} 
We shall give only the proof for the case $\IF>x_2>0>-x_1>-\IF$.
The other cases can be established by similar arguments.
Since our approach is \tb{of asymptotic nature},
we assume in the following that $v$ is sufficiently large \tb{so} that $x_i(v)>0, i=1,2$. Let $S>2\lambda$ be any positive constant. With $q=q(v)$ defined in \eqref{eq:qq} we denote
$$
\Del_k=\LT[kS q,(k+1)S q\RT],\  k\in \mathbb{Z}, \ \ \text{and}\ \
 N_i(v)=\LT\lfloor S^{-1} x_i(v)q^{-1} v^{-1}\RT\rfloor,\ i=1,2,
$$
 where $\lfloor\cdot\rfloor$ is the ceiling function.
For any small $\vn_0>0$, set $\lamp=\lambda+{\vn_0}$ and $\lamm=\max(0,\lambda-{\vn_0})$.
It follows by Bonferroni's inequality \tb{that}
\BQN \label{eq:WPP}
\sum_{k=-N_1(v)-1}^{N_2(v)+1} Q_k^+(v)\ge \pk{\sup_{t\in\Theta_{x_1,x_2}(v)}\inf_{s\in[0,\lambda(v)  q]}Z(t+s)>v} \ge  \sum_{k=-N_1(v)}^{N_2(v)} Q_k^-(v) -\Sigma_1(v)
\EQN
for large enough $u$, where
\BQNY
Q_k^+(v)&=&\pk{\sup_{t\in\Del_k}\inf_{s\in\LT[0, \lamm  q \RT]} Z(t_0+t+s) >v},\ \ k\in\mathbb{Z},\\
Q_k^-(v)&=&\pk{\sup_{t\in\Del_k}\inf_{s\in\LT[0, \lamp  q \RT]} Z(t_0+t+s) >v},\ \ k\in\mathbb{Z},\\
\Sigma_1(v)&=& \underset{ -N_1(v)\le k<l\le N_2(v)}{\sum}\pk{\sup_{t\in\Del_k}\inf_{s\in\LT[0, \lamp  q \RT]} Z(t_0+t+s) >v,\sup_{t\in\Del_l}\inf_{s\in\LT[0, \lamp  q \RT]}Z(t_0+t+s) >v}.
\EQNY
Next, we shall derive upper bounds for $Q_k^+(v)$ and lower bounds for $Q_k^-(v)$. First, note that 
\BQNY
 Q_k^+(v)&\le&\pk{\sup_{t\in\Del_k}\inf_{s\in\LT[0, \lamm  q \RT]} \overline Z(t_0+t+s) >\frac{v}{\sigma_Z^+(k,v)}}  \\
Q_k^-(v)&\ge& \pk{\sup_{t\in\Del_k}\inf_{s\in\LT[0, \lamp  q \RT]} \overline Z(t_0+t+s)>\frac{v} {\sigma_Z^-(k,v)}},
\EQNY
where $\overline Z(t):=Z(t)/\sigma_Z(t), t\ge0$ and
$$
 \sigma_Z^- (k,v)=\inf_{t\in\Del_k}\inf_{s\in\LT[0, \lamp  q \RT]}\sigma_Z(t_0+t+s),\ \
 \sigma_Z^+(k,v)=\sup_{t\in\Del_k}\sup_{s\in\LT[0, \lamm  q \RT]}\sigma_Z(t_0+t+s).
$$
Furthermore, since
\BQN\label{eq:sig1}
 \sigma_{Z}(t)=A-\frac{A^2 B}{2}(t-t_0)^2\oo,\ \ \ t\to t_0,
 \EQN
for any  small $\vn_1>0$ there exists $v_0$ such that for any $v>v_0$  \Ea{(set below $B^{\pm}=B(1 \pm \ve_1)$)} 
\BQNY
\frac{1}{\sigma_Z^-(k,v)}\le \frac{1}{A}+\frac{\Ea{ B^+}}{2}\LT(((k+1)S+\lamp) q\RT)^2,\quad
\frac{1}{\sigma_Z^+(k,v)}\ge \frac{1}{A}+\frac{B^{-}}{2}\LT(kS  q\RT)^2 
\EQNY
hold  for $ k=0,\cdots,N_2(v)+1$, and also
\BQNY
\frac{1}{\sigma_Z^-(k,v)}\le \frac{1}{A}+\frac{\Ea{ B^+}}{2}\LT(kS  q\RT)^2,\quad
\frac{1}{\sigma_Z^+(k,v)}\ge \frac{1}{A}+\frac{\Ea{ B^-}}{2}\LT(((k+1)S+\lamm) q\RT)^2 
\EQNY
hold  for  $k=-N_1(v)-1,\cdots,-1$.
Moreover, for any $k=-N_1(v)-1,\cdots,N_2(v)+1$,
set \tb{$\overline Z_{k,v}(t,s)=\overline Z(t_0+kS q+t q+s q), (t,s)\in[0,S]\times[0,\lamp]$}.
It follows from \eqref{eq:locstaXH} that, for the correlation function $r_{\overline Z_{k,v}}(\cdot,\cdot,\cdot,\cdot)$  of $\overline Z_{k,v}$
\BQN
\lim_{v\to\IF}  2v^2(1-r_{\overline Z_{k,v}}(t,s,t',s'))= Q\abs{s+t-s'-t'}^{\alpha}
\EQN
holds for any $(t,s),(t',s')\in[0,S]\times[0,\lamp]$. Furthermore, for sufficiently large $v$
\BQNY
2v^2(1-r_{\overline Z_{k,v}}(t,s,t',s'))&\le& G_0\frac{K^2(q \abs{s+t-s'-t'})}{K^2(q)},
\EQNY
for all $(t,s),(t',s')\in[0,S]\times[0,\lamp]$, with some positive constant $G_0$. Set $S_{\max}=\max\{\abs{s+t-s'-t'}: (t,s),(t',s')\in[0,S]\times[0,\lamp]\}$. Using Potter bounds (cf. \cite{EKM97}), for any small $\delta>0$ we have, when $v$ is sufficiently large
\BQNY
\frac{K^2(q \abs{s+t-s'-t'})}{K^2(q)}&\le &G_1 \max\LT(S_{\max}^{\alpha-\delta},S_{\max}^{\alpha+\delta} \RT) \LT(\frac{\abs{s+t-s'-t'}}{S_{\max}}\RT)^{\alpha-\delta}\\
&\le& G_2 (\abs{t-t'}^{\alpha-\delta}+\abs{s-s'}^{\alpha-\delta})
\EQNY
holds uniformly with respect to $(t,s),(t',s')\in[0,S]\times[0,\lamp]$, where $G_1,G_2$ are two positive constants.
Hence, by an application of \nelem{LemGP}, where we set
$$f(v)=\frac{v}{A},\ \ \theta(v)=\LT(\frac{1}{A}+\frac{ \Ea{ B^+}}{2}\LT(((k+1)S+\lamp) q\RT)^2\RT) v, \ \ D=\frac{Q}{2 A^2},$$
we obtain, for any $k=0,\cdots,N_2(v)+1$
\BQNY
Q_k(v)\ge  \H_\alpha(D_0 S, D_0 \lamp)\frac{1}{\sqrt{2\pi}\theta(v)}\exp\LT(-\frac{(\theta(v))^2}{2}\RT)\oo, \ \ u\to\IF, 
\EQNY
where $D_0=D^{\frac{1}{\alpha}}=2^{-\frac{1}{\alpha}} A^{-\frac{2}{\alpha}} Q^{\frac{1}{\alpha}}.$
Therefore, as $v\to\IF$ \Ea{(set below $\zeta(v)=v^{-2}q^{-1}\exp(-\frac{ v^2}{2 A^2}) $)}
\BQN\label{eq:low1}
\sum_{k=0}^{N_2(v)}Q_k(v)&\ge&  \H_\alpha(D_0 S, D_0 \lamp) \frac{ A }{\sqrt{2\pi}v} \sum_{k=0}^{N_2(v)} \exp\LT(-\frac{v^2 \LT(\frac{1}{A}+\frac{ \Ea{ B^+}}{2}\LT(((k+1)S+\lamp) q\RT)^2\RT)^2  }{2}\RT)\oo\nonumber\\
&=& \frac{1}{S}\H_\alpha(D_0 S, D_0 \lamp) \frac{ A}{\sqrt{2\pi}}\Ea{\zeta(v)} \int_0^{x_2} \exp\LT(- \frac{\Ea{ B^+}}{2 A}x^2   \RT)dx \oo,
\EQN
where we \tb{used  that} 
$\lim_{v\to\IF} v q=\lim_{v\to\IF} v \ovSK\LT(\frac{1}{v}\RT)=0$ and
$ \lim_{v\to\IF}x_2(v)v^{- {1}/{2}}=0$.

Similarly,  as $v\to\IF$
\BQN\label{eq:low2}
\sum_{k=-N_1(v)}^{-1}Q_k(v)&\ge&
\frac{1}{S}\H_\alpha(D_0 S, D_0 \lamp) \frac{ A }{\sqrt{2\pi}}\Ea{\zeta(v)}  \int_{-x_1}^0 \exp\LT(- \frac{\Ea{B^+}}{2 A}x^2   \RT)dx \oo.
\EQN
Furthermore, with the same arguments as above for any $S_1>2\lambda$
\BQN\label{eq:upper}
\sum_{k=-N_1(v)-1}^{N_2(v)+1}Q_k(v)&\le&
\frac{1}{S_1}\H_\alpha(D_0 S_1, D_0 \lamm) \frac{ A }{\sqrt{2\pi}} \Ea{\zeta(v)}
\int_{-x_1}^{x_2} \exp\LT(- \frac{\Ea{B^-}}{2 A}x^2   \RT)dx \oo.
\EQN
Consequently, \eqref{eq:WPP} and (\ref{eq:low1}-\ref{eq:upper}) imply \Ea{(set $\bar \zeta(v):= D_0A^{\frac32}\zeta(v)/ \sqrt{B^+}$)}
\BQN\label{eq:HH}
\lefteqn{\frac{1}{D_0 S_1}\H_\alpha(D_0 S_1, D_0 \lamm) \LT(\Phi\LT(\LT(\frac{\Ea{B^-}}{A}\RT)^{\frac12}x_2\RT)-\Phi\LT(-\LT(\frac{\Ea{B^-}}{A}\RT)^{\frac12}x_1\RT)\RT)}\nonumber\\
&  \ge &\limsup_{v\to\IF}  
\pk{\sup_{t\in\Theta_{x_1,x_2}(v)}\inf_{s\in[0,\lamm q]}Z(t+s)>v}/\lp{\bar \zeta(v)}
\nonumber\\
&  \ge &\limsup_{v\to\IF}  
\pk{\sup_{t\in\Theta_{x_1,x_2}(v)}\inf_{s\in[0,\lambda(v) q]}Z(t+s)>v}/\lp{\bar \zeta(v)}
\nonumber\\
&  \ge &\liminf_{v\to\IF}   
\pk{\sup_{t\in\Theta_{x_1,x_2}(v)}\inf_{s\in[0,\lambda(v) q]}Z(t+s)>v}/\lp{\bar \zeta(v)}
\nonumber\\
&  \ge &\liminf_{v\to\IF}   
\pk{\sup_{t\in\Theta_{x_1,x_2}(v)}\inf_{s\in[0,\lamp q]}Z(t+s)>v}/\lp{\bar \zeta(v)}
\\
&\ge  &\frac{1}{D_0 S}\H_\alpha(D_0 S, D_0 \lamp) \LT(\Phi\LT(\LT(\frac{\Ea{B^+}}{A}\RT)^{\frac12}x_2\RT)-\Phi\LT(-\LT(\frac{\Ea{B^+}}{A}\RT)^{\frac12}x_1\RT)\RT)
 -\limsup_{v\to\IF}   
\Sigma_1(v)/\lp{\bar \zeta(v)}
.\nonumber
\EQN
Moreover, since
\BQNY
\Sigma_1(v)&\le& \underset{ -N_1(v)\le k<l\le N_2(v)}{\sum}\pk{\sup_{t\in\Del_k}  {Z}(t_0+t) >v,\sup_{t\in\Del_l} {Z}(t_0+t) >v}
\EQNY
similar arguments as in the proof of Eqs. (31) and (32)  in \cite{HJ13} imply
\BQN\label{eq:sig0}
\lim_{S\to\IF} \limsup_{v\to\IF}   \Sigma_1(v)/\lp{\bar \zeta(v)} 
=0.
\EQN
Let us assume for the moment that
\BQN\label{eq:LowP}
\limsup_{S\to\IF}\frac{1}{S}\H_\alpha(S, D_0 \lambda)>0.
\EQN
Letting first $\vn_0, \vn_1\to 0$ and then $S,S_1\to\IF$ we get from \eqref{eq:HH} \Ea{and the definition of $\H_\alpha$} that
\BQNY
\IF>\H_\alpha\ge \liminf_{S\to\IF}\frac{1}{S}\H_\alpha(S, D_0 \lambda)\ge\limsup_{S\to\IF}\frac{1}{S}\H_\alpha(S, D_0 \lambda)>0.
\EQNY
Further, in view of \eqref{eq:HPA1} and \eqref{eq:HPA2} we have
$$
\pk{\sup_{t\ge0} Z(t)>v }=D_0A^{\frac32}B^{-\frac12}\H_\alpha  \Ea{\zeta(v)} 
\oo,\ \ \text{as}\ v\to\IF.
$$
Therefore,  the claim of \netheo{ThmX}  follows with $\GH_\alpha (\lambda)\in(0,\IF)$.\\
 Next, we prove \eqref{eq:LowP}. Define
$$
E_v=\bigcup_k \Biggl(\Del_{2k}\cap \Theta_{x_1,x_2}(v)\Biggr), \quad
N^*(v)=\sharp\{k\in\mathbb{Z}: \Del_{2k}\cap \Theta_{x_1,x_2}(v)\neq \emptyset\}.
$$
For any $v$ positive
\BQN\label{eq:Ev}
\pk{\sup_{t\in\Theta_{x_1,x_2}(v)}\inf_{s\in[0,\lamp  q]}Z(t,s)>v}\ge \pk{\sup_{t\in E_v}\inf_{s\in[0,\lamp  q]}Z(t,s)>v}.
\EQN
Using Bonferroni's inequality and the same arguments as in the derivation of \eqref{eq:low1} we conclude that
\BQN
 \pk{\sup_{t\in E_v}\inf_{s\in[0,\lamp  q]}Z(t,s)>v}
&\ge &
\frac{1}{2S}\H_\alpha(D_0 S, D_0 \lamp) \frac{ A }{\sqrt{2\pi}}\Ea{\zeta(v)}
\int_{-x_1}^{x_2} \exp\LT(- \frac{\Ea{B^+}}{2 A}x^2   \RT)dx- \Sigma_2(v),
\EQN
where
\BQNY
\Sigma_2(v)&=& \underset{k,l\in N^*(v), k>l}{\sum} \pk{\sup_{t\in\Del_{2k}}\inf_{s\in\LT[0, \lamp  q \RT]} Z(t_0+t+s) >v,\sup_{t\in\Del_{2l}}\inf_{s\in\LT[0, \lamp  q \RT]}Z(t_0+t+s) >v}\\
&\le&\underset{k,l\in N^*(v), k>l}{\sum} \pk{\sup_{t\in\Del_{2k}}  Z(t_0+t) >v,\sup_{t\in\Del_{2l}}Z(t_0+t) >v}.
\EQNY
Similar arguments as in the proof of Eq. (32)  in \cite{HJ13} show that
\BQN\label{eq:sig00}
 \limsup_{v\to\IF}  
 \Sigma_1(v) /\lp{\bar \zeta(v)}
 \le G_3 S \sum_{k\ge1}\exp\LT(-G_4(k S)^{\alpha}\RT)
\EQN
for some positive constants $G_3, G_4$.
Therefore, combining \eqref{eq:HH}, (\ref{eq:Ev}-\ref{eq:sig00}) we conclude that
$$
\liminf_{S_1\to\IF}\frac{1}{ S_1}\H_\alpha( S_1, D_0 \lambda)\ge \frac{1}{ S}\LT(\frac{1}{2D_0}\H_\alpha( D_0S, D_0 \lambda)-G_5 S^2 \sum_{k\ge1}\exp\LT(-G_4(k S)^{\alpha}\RT)\RT),
$$
with some positive constant $G_5$.
Since $\H_\alpha( D_0S, D_0 \lambda)$ is positive and increasing as $S$ increases, then for $S$ sufficiently large the right hand side in the last formula is strictly positive, implying thus \eqref{eq:LowP}. This completes the proof. \QED

\bigskip

{\bf Acknowledgement}: The authors kindly acknowledge partial support by
the Swiss National Science Foundation Grant 200021-140633/1, and the project RARE -318984 (a Marie Curie IRSES FP7 Fellowship).
The first author also acknowledges
partial support by NCN Grant No 2013/09/B/ST1/01778 (2014-2016).

\bibliographystyle{plain}

 \bibliography{gausbibruinABCD}
\end{document}